\documentclass[fleqn]{mat01}
\usepackage{times,mathtimy,amssymb}
\begin{document}

\setcounter{page}{55}
\firstpage{55}

\font\xx=msam5 at 10pt
\def\ab{\mbox{\xx{\char'03}}}

\newtheorem{theo}{Theorem}
\renewcommand\thetheo{\arabic{section}.\arabic{theo}}
\newtheorem{theor}[theo]{\bf Theorem}
\newtheorem{lem}[theo]{Lemma}
\newtheorem{propo}[theo]{\rm PROPOSITION}
\newtheorem{rema}[theo]{Remark}
\newtheorem{defn}[theo]{\rm DEFINITION}
\newtheorem{exam}{Example}
\newtheorem{pol}[theo]{Proof of Lemma}
\newtheorem{coro}[theo]{\rm COROLLARY}
\newtheorem{claim}[theo]{Claim}
\newtheorem{conjecture}[theo]{Conjecture}

\def\h{{\cal H}}
\def\n{{\cal N}}
\def\m{{\cal M}}
\def\u{{\cal U}}
\def\a{{\cal A}}
\def\l{{\cal L}}
\def\b{{\cal B}}
\def\s{{\cal S}}
\def\t{{\cal T}}
\def\r{{\cal R}}
\def\p{{\cal P}}
\def\bh{{\cal B({\cal H})}}
\def\fh{{\cal F(H)}}
\def\ch{{\cal C}_{1}(\h)}
\def\alg{\mbox{\rm Alg\,}}
\def\lat{\mbox{\rm Lat\,}}
\def\hom{\mbox{\rm Hom\,}}
\def\homp{\mbox{\rm Hom}(\lat\a, \p)}
\def\ref{\mbox{\rm Ref\ }}
\def\uir{{\cal U}_{\phi}^{r}}
\def\uil{{\cal U}_{\phi}^{l}}
\def\ui{{\cal U_{\phi}}}
\def\ut{{\cal U}_{\tau}}
\def\ri{{\cal R}_{\phi}}
\def\rt{{\cal R}_{\tau}}
\def\ps{\phi_{\sim}}
\def\pss{(\phi_{\sim})_{\sim}}
\def\pis{\psi_{\sim}}
\def\piss{(\psi_{\sim})_{\sim}}
\def\ts{\tau_{\sim}}
\def\tss{(\tau_{\sim})_{\sim}}
\def\bo{\overline{\otimes}}
\def\pf{{\bf Proof.}\,}
\def\ox{$\qquad\Box$}
\def\q{\quad}
\def\qq{\qquad}
\def\sub{\subseteq}
\def\si{\sigma}
\def\pa{\parallel}

\font\vvv=cmbsy10 at 10.4pt
\def\aca{\mbox{\vvv{\char'101}}}

\font\uuu=cmbsy10 at 13.5pt
\def\acac{\mbox{\uuu{\char'101}}}

\title{Rank-one operators in reflexive one-sided $\acac$-submodules}

\markboth{Dong Zhe}{Rank-one operators in reflexive one-sided $\a$-submodules}

\author{DONG ZHE}

\address{Department of Mathematics, Zhejiang University,
Hangzhou 310 027, People's~Republic of China\\
\noindent E-mail: dongzhe@zju.edu.cn}

\volume{114}

\mon{February}

\parts{1}

\Date{MS received 27 May 2003; revised 16 October 2003}

\begin{abstract}
In this paper, we first characterize reflexive one-sided $\a$-submodules
$\u$ of a unital operator algebra $\a$ in $\bh$ completely. Furthermore
we investigate the invariant subspace lattice $\lat\r$ and the reflexive
hull $\ref\r$, where $\r$ is the submodule generated by rank-one
operators in $\u$; in particular, if $\l$ is a subspace lattice, we
obtain when the rank-one algebra $\r$ of $\alg\l$ is big enough to
determined $\alg\l$ in the following senses: $\alg\l=\alg\lat\r$ and
$\alg\l=\ref\r$.
\end{abstract}

\keyword{Reflexive one-sided $\a$-submodule; rank-one operator.}

\maketitle

\section{Introduction}

Let $\h$ be a complex Hilbert space, $\bh$ the algebra of all bounded
linear operators on $\h$ and $\p$ the complete lattice of all orthogonal
projections in $\bh$. Suppose that $\a$ is a unital operator algebra in
$\bh$ and $\phi$ is an order homomorphism of $\lat\a$ into itself (i.e.
$E\leq F$ implies $\phi(E)\leq\phi(F)$), where $\lat\a$ is the complete
lattice of all invariant projections for $\a$. Then the set
$\u=\{T\in\bh: TE\sub\phi(E)\;\mbox{for all}\; E\in\lat\a\}$ is clearly
a weakly closed two-sided $\a$-submodule of $\bh$. 

It became apparent that many interesting classes of non-self adjoint
operator algebras arise as just such a module. Erdos and Power
in \cite{3} proved that any weakly closed $\a$-submodule of $\bh$ for a nest
algebra $\a$ is of the above form. In \cite{4}, Han Deguang proved that this
is also true for any reflexive algebra $\a$, which is $\si$-weakly
generated by rank-one operators in itself. The purpose of this paper is
to show that any reflexive right $\a$-submodule and $*$-reflexive left
$\a$-submodule of a unital operator algebra $\a$ are determined by order
homomorphisms from $\lat\a$ into $\p$. As a corollary, we obtain the
complete characterization of all $\si$-weakly closed one-sided
$\a$-submodules, where $\a$ is $\si$-weakly generated by rank-one
operators in itself or, in particular, $\a$ is a nest algebra.

In \cite{2}, Erdos showed that if $\lat\a$ is a nest then the set of
finite sums of rank-one operators in $\a$ is $\si$-weakly dense in
$\a$. In \cite{9}, Longstaff asked whether the same conclusion holds
for the more general case of completely distributive lattices, and
showed that, in the opposite direction, complete distributivity is
a necessary condition for this. Subsequently, Lambrou \cite{6}
showed that complete distributivity of the invariant subspace
lattices implies a condition somewhat weaker than the strong
density. Laurie and Longstaff \cite{7} proved that the answer is
affirmative if additional requirement of commutativity is
imposed on the invariant subspace lattice. In \S3, we will
consider when the rank-one subalgebra $\r$ of $\alg\l$ determines
$\alg\l$ in senses other than the $\si$-weak density.

Which subspace lattices $\l$ are determined by the rank-one
subalgebra $\r$ of $\alg\l$ in the sense that $\l=\lat\r$? This
question was answered by Longstaff in (\cite{8}, Proposition~3.2). A
sufficient but not necessary condition (\cite{8}, Corollary~3.2.1) was
given and it is shown in \cite{8} that this condition is strictly
weaker than complete distributivity. In \S3, we investigate
the invariant subspace lattice of the rank-one submodule of $\u$.
As an application, we derive the sufficient and necessary
condition obtained by Longstaff in \cite{8} in order that
$\l=\lat\r$. As another application, we also obtain an equivalent
condition for which $\alg\l=\alg\lat\r$.

In \S3, we also study when the rank-one submodule $\r$ of a reflexive
one-sided $\a$-submodule $\u$ is big enough to determine $\u$ in the
sense that $\ref\r=\u$, where $\ref\r=\{T\in\bh: Tx\in [\r x]\;\mbox{for
all}\; x\in\h\}$ is the reflexive hull of $\r$. An equivalent condition
for $\ref\r=\u$ is given by means of order homomorphisms from $\lat\a$
into $\p$.

The terminology and notation of this paper concerning reflexive
subspaces may be found in \cite{5}. In what follows, we always assume
that $\a$ is a unital operator algebra in $\bh$. Set
\begin{equation*}
\hom(\lat\a, \p)=\{\phi: \phi\;\mbox{is an order homomorphism
from}\;\lat\a\; \mbox{into}\; \p\}.
\end{equation*}

$\left.\right.$\vspace{-1.6pc}

\noindent Given $\phi$ in $\homp$, a right $\a$-submodule is associated
which is given by
\begin{equation*}
\u_{\phi}^{r}=\{T\in\bh: TE\sub\phi(E), \forall E\in\lat\a\};
\end{equation*}
and a left $\a$-submodule which is given by
\begin{equation*}
\u_{\phi}^{l}=\{T\in\bh: T\phi(E)\sub E, \forall E\in\lat\a\}.
\end{equation*}
Clearly they are weakly closed. We say that $\uir$ (and $\uil$) are the
right(left) $\a$-submodule determined by $\phi$ respectively. To each
$\phi$ in $\homp$ there is naturally associated $\ps$ in $\homp$ given
by
\begin{equation*}
\ps(E)=\vee\{F\in\lat\a: \phi(F)\not\geq E\},\q\forall E\in\lat\a
\end{equation*}
(with the convention that $\ps(0)=0$). Observe that $\homp$ has a
natural partial ordering given by $\phi\leq\psi$ if and only if
$\phi(E)\leq\psi(E)$ for any $E\in\lat\a$. It follows that
$\phi\leq\psi$ implies $\ps\geq\pis$.

\section{Basic properties of one-sided $\aca$-submodules}

A subspace $\s$ of $\bh$ is said to be $*$-reflexive, if $\s^{*}$ is
reflexive.

\begin{theor}[\!] Suppose that $\a$ is a unital operator
algebra in $\bh$ and $\u$ is a subspace of $\bh$. Then

\begin{enumerate}
\renewcommand\labelenumi{\rm (\arabic{enumi})}
\item $\u$ is  a reflexive right $\a$-submodule if and only if
there exists $\phi\in\homp$ such that
$\u=\{T\in\bh: TE\sub\phi(E), \forall E\in\lat\a\};$

\item $\u$ is a $*$-reflexive left $\a$-submodule if and only if
there exists $\psi\in\homp$ such that $\u=\{T\in\bh: T\psi(E)\sub
E, \forall E\in\lat\a\}.$\vspace{-.7pc}
\end{enumerate}
\end{theor}

\begin{proof}
(1) {\it Sufficiency.}\ \ Clearly $\u$ is a right $\a$-submodule, so
we only need to prove that $\u$ is reflexive. Suppose that $T\in\bh$ and
$Tx\in [\u x]$ for any $x\in\h$. Thus for any $E\in\lat\a$,
\begin{equation*}
TE\sub [\u E]=[\phi(E)\u E]=\phi(E)[\u E]\sub\phi(E).
\end{equation*}
So $T\in\u$ and it shows that $\u$ is reflexive.

\newpage

{\it Necessity.}\ \ For any $E\in\lat\a$, let $\phi(E)=[\u E]$. Clearly
$\phi$ is an order homomorphism in $\homp$. Set
\begin{equation*}
\uir=\{T\in\bh: TE\sub\phi(E), \forall E\in\lat\a\}.
\end{equation*}
It is obvious that $\u\sub\uir$. Conversely, let $T\in\uir$. For any
$x\in\h$, denote by $E$ the orthogonal projection onto $[\a x]$. Then
$E\in\lat\a$, $x\in E$ and
\begin{equation*}
Tx\in TE\sub\phi(E)=[\u E]=[\u[\a x]]=[\u x]
\end{equation*}
since $\u$ is a right $\a$-submodule. From the reflexivity of $\u$, it
follows that $T\in\u$. Accordingly, $\uir\sub\u$ and $\u=\uir$.

(2) {\it Sufficiency.}\ \ Suppose that there exists $\psi\in\homp$ such
that
\begin{equation*}
\u=\{T\in\bh: T\psi(E)\sub E, \forall E\in\lat\a\}.
\end{equation*}
Define $\phi: \lat\a^{*}=(\lat\a)^{\bot}\rightarrow\p$ by
\begin{equation*}
\phi(E^{\bot})=I-\psi(E),\q\forall E^{\bot}\in\lat\a^{*}=(\lat\a)^{\bot}.
\end{equation*}
Certainly $\phi\in\hom(\lat\a^{*}, \p)$. Thus
\begin{align*}
\u^{*} &= \{T^{*}\in\bh: T\psi(E)\sub E, \forall E\in\lat\a\}\\
&=\{T^{*}\in\bh: T^{*}E^{\bot}\sub\psi(E)^{\bot}=\phi(E^{\bot}), \forall
E\in\lat\a\}\\
&= \{S\in\bh: SE^{\bot}\sub\phi(E^{\bot}), \forall E^{\bot}\in\lat\a^{*}
= (\lat\a)^{\bot}\}.
\end{align*}
It follows from (1) that $\u^{*}$ is a reflexive right
$\a^{*}$-submodule, and $\u$ is a $*$-reflexive left $\a$-submodule.

{\it Necessity.}\ \ Suppose that $\u$ is a $*$-reflexive left
$\a$-submodule. Thus $\u^{*}$ is a reflexive right $\a^{*}$-submodule,
it follows from (1) that there exists $\phi\in\hom(\lat\a^{*}, \p)$ such
that
\begin{equation*}
\u^{*}=\{T\in\bh: TE^{\bot}\sub\phi(E^{\bot}), \;\forall
E^{\bot}\in\lat\a^{*}\}.
\end{equation*}
Define $\psi: \lat\a\rightarrow\p$ by $\psi(E)=I-\phi(E^{\bot})$.
Clearly $\psi\in\homp$ and
\begin{align*}
\u &= \{T^{*}\in\bh: T^{*}\phi(E^{\bot})^{\bot}\sub E, \;\forall
E\in\lat\a\}\\
&= \{S\in\bh: S\psi(E)\sub E, \;\forall E\in\lat\a\}.
\end{align*}

$\left.\right.$\vspace{-1.5pc}

\hfill \ab
\end{proof}

From the proof of Theorem~2.1, we know that if $\u$ is a reflexive right
$\a$-submodule then $\u=\{T\in\bh: TE\sub\tau_{r}(E), \forall
E\in\lat\a\}$, where $\tau_{r}(E)=[\u E]$; if $\u$ is a $*$-reflexive
left $\a$-submodule then $\u=\{T\in\bh: T\tau_{l}(E)\sub E, \forall
E\in\lat\a\}$, where $\tau_{l}(E)=I-[\u^{*}E^{\bot}]$.

\begin{coro}$\left.\right.$\vspace{.5pc}

\noindent If $\a$ is a unital $\si$-weakly closed algebra which is
$\si$-weakly generated by rank-one operators in $\a${\rm ,} then every
$\si$-weakly closed right or left $\a$-submodule has the form given in
Theorem~{\rm 2.1(1)} or {\rm (2),} respectively.
\end{coro}

\begin{proof} By virtue of (\cite{5}, Theorem~2.2), every $\si$-weakly
closed right or left $\a$-submodule is reflexive. So the result is true
for $\si$-weakly closed right $\a$-submodule by Theorem~2.1(1). Now for
any $\si$-weakly closed left $\a$-submodule $\u$, since the adjoint
operation is continuous in the $\si$-weak topology, $\u^{*}$ is a
$\si$-weakly closed right $\a^{*}$-submodule and $\a^{*}$ is
$\si$-weakly generated by rank-one operators in $\a^{*}$. Therefore it
follows from (\cite{5}, Theorem~2.2) that $\u^{*}$ is reflexive and $\u$
is $*$-reflexive. Thus $\u$ has the form in Theorem~2.1(2).\hfill \ab
\end{proof}

\begin{coro}$\left.\right.$\vspace{.5pc}

\noindent Suppose that $\l$ is a commutative and completely distributive
subspace lattice{\rm ,} or specially{\rm ,} a nest. Then every
$\si$-weakly closed right or left $\alg\l$-submodule is of the form
given in Theorem~{\rm 2.1(1)} or {\rm (2),} respectively.
\end{coro}

\begin{proof}
This follows from Corollary~2.2 and (\cite{7}, Theorem~3).\hfill\ab
\end{proof}

\begin{coro}$\left.\right.$\vspace{.5pc}

\noindent Suppose that $\a$ is a unital algebra in $\bh$.
\begin{enumerate}
\renewcommand\labelenumi{\rm (\arabic{enumi})}
\item Let $\u$ be as in {\rm (1)} of Theorem~{\rm 2.1.} Then $\u$ is a
right ideal if and only if $\tau_{r}(E)\leq E$ for every
$E\in\lat\a${\rm ,} where $\tau_{r}(E)=[\u E]${\rm ;} 

\item Let $\u$ be as in {\rm (2)} of Theorem~{\rm 2.1.} Then $\u$ is a
left ideal if and only if $\tau_{l}(E)\geq E$ for any $E\in\lat\a${\rm
,} where $\tau_{l}(E)=I-[\u^{*}E^{\bot}]$.\vspace{-.5pc}
\end{enumerate}
\end{coro}

\begin{proof}$\left.\right.$

\begin{enumerate}
\renewcommand\labelenumi{\rm (\arabic{enumi})}
\item Obvious.

\item Let $\u$ be a left ideal of $\a$. Thus $\u^{*}$ is a right ideal
of $\a^{*}$, it follows from (1) that $[\u^{*}E^{\bot}]\leq E^{\bot}$
for any $E^{\bot}\in\lat\a^{*}=(\lat\a)^{\bot}$. This deduces that
$\tau_{l}(E)=I-[\u^{*}E^{\bot}]\geq E$ for any $E\in\lat\a$. The
converse implication can be proved similarly.

\hfill \ab
\end{enumerate}\vspace{-1.5pc}
\end{proof}

\begin{propo}$\left.\right.$\vspace{.5pc}

\noindent Suppose that $\a$ is a unital algebra in $\bh$.

\begin{enumerate}
\renewcommand\labelenumi{\rm (\arabic{enumi})}
\item Let $\u$ be a reflexive right $\a$-submodule. Then $P\in\lat\u$ if
and only if there exists $E\in\lat\a$ such that $\tau_{r}(E)\leq P\leq
E${\rm ;}

\item Let $\u$ be a $*$-reflexive left $\a$-submodule. Then $P\in\lat\u$
if and only if there exists $E\in\lat\a$ such that $E\leq
P\leq\tau_{l}(E)$.\vspace{-.5pc}
\end{enumerate}
\end{propo}

\begin{proof}$\left.\right.$

\begin{enumerate}
\renewcommand\labelenumi{\rm (\arabic{enumi})}
\item From the proof of Theorem~2.1, $\u=\{T\in\bh: TE\sub\tau_{r}(E),
E\in\lat\a\}$. If $\tau_{r}(E)\leq P\leq E$ for some $E\in\lat\a$ and
$T\in\u$, then
\begin{equation*}
\hskip -1.25pc TP=TEP=\tau_{r}(E)TEP=P\tau_{r}(E)TEP=PTP.
\end{equation*}
So $P\in\lat\u$.

\hskip 1pc Conversely, if $P\in\lat\u$, let $E=[\a P]$. Then
$E\in\lat\a$, $E\geq P$ and 
\begin{equation*}
\hskip -1.25pc \tau_{r}(E)=[\u E]=[\u[\a P]]\sub [\u P]\sub P
\end{equation*}
since $\u$ is a right $\a$-module. Thus $\tau_{r}(E)\leq P\leq E$.

\item Follows from (1) and a simple calculation.\hfill \ab
\end{enumerate}\vspace{-2.5pc}
\end{proof}

\newpage

For non-zero vectors $x, y\in\h$, the rank-one operator $x\otimes y$ is
defined by the equation
\begin{equation*}
(x\otimes y)z = \langle z, y\rangle x,\q\forall z\in\h.
\end{equation*}

\begin{lem} Suppose that $\a$ is a unital algebra in $\bh$.

\begin{enumerate}
\renewcommand\labelenumi{\rm (\arabic{enumi})}
\item Let $\u^{r}_{\phi}$ be the reflexive right $\a$-submodule
determined by $\phi$ in $\homp$. Then a rank-one operator $x\otimes
y\in\u^{r}_{\phi}$ if and only if for some $E\in\lat\a$, $x\in E$ and
$y\in\ps(E)^{\bot}${\rm ,} where $\ps(E)=\vee\{F\in\lat\a: \phi(F)\not\geq
E\}$.

\item Let $\u^{l}_{\phi}$ be the $*$-reflexive left $\a$-submodule
determined by $\phi$ in $\homp$. Then a rank-one operator $x\otimes
y\in\u^{l}_{\phi}$ if and only if for some $E\in\lat\a${\rm ,}
$x\in\wedge\{F\in\lat\a: \phi(F)\not\leq E\} $ and $y\in E^{\bot}$.
\end{enumerate}
\end{lem}

\begin{proof}$\left.\right.$

\begin{enumerate}
\renewcommand\labelenumi{\rm (\arabic{enumi})}
\item Suppose that there exists $E\in\lat\a$ such that $x\in E$ and
$y\in\ps(E)^{\bot}$. For any $F\in\lat\a$, if $\phi(F)\geq E$, then
\begin{equation*}
\hskip -1.25pc (x\otimes y)F=E(x\otimes y)\ps(E)^{\bot}F\sub E\sub\phi(F);
\end{equation*}
if $\phi(F)\not\geq E$, it follows from the definition of $\ps(E)$ that
$F\leq\ps(E)$. Thus
\begin{equation*}
\hskip -1.25pc (x\otimes y)F = E(x\otimes y)\ps(E)^{\bot}F = 0\sub\phi(F).
\end{equation*}
Accordingly, $x\otimes y\in\u^{r}_{\phi}$.

\hskip 1pc Conversely, if $x\otimes y\in\u^{r}_{\phi}$. Let
\begin{equation*}
\hskip -1.25pc E = \wedge\{F\in\lat\a: Fx=x\}.
\end{equation*}
Naturally, $E\in\lat\a$ and $x\in E$. For any $F\in\lat\a$ and
$\phi(F)\not\geq E$, it follows from the definition of $E$ that
$\phi(F)x\not=x$. Since $x\otimes y\in\u^{r}_{\phi}$, we have
\begin{equation*}
\hskip -1.25pc (x\otimes y)Fy = \phi(F)(x\otimes y)Fy
\end{equation*}
and
\begin{equation*}
\hskip -1.25pc \|Fy\|^{2}x = \|Fy\|^{2}\phi(F)x.
\end{equation*}
So $Fy=0$. From the definition of $\ps(E)$, it follows that $\ps(E)y=0$
and $y\in\ps(E)^{\bot}$.

\item By hypothesis, $\u^{l}_{\phi}=\{T\in\bh: T\phi(E)\sub E, \forall
E\in\lat\a\}$. Define $\psi: \lat\a^{*}\rightarrow\p$ by
$\psi(E^{\bot})=I-\phi(E)$. Thus $\psi\in\hom(\lat\a^{*}, \p)$ and
\begin{align*}
\hskip -1.25pc (\u^{l}_{\phi})^{*} &= \{T^{*}\in\bh: T\phi(E)\sub E, \forall E\in\lat\a\}\\
\hskip -1.25pc &= \{T^{*}\in\bh: T^{*}E^{\bot}\sub\phi(E)^{\bot}, \forall E\in\lat\a\}\\
\hskip -1.25pc &= \{S\in\bh: SE^{\bot}\sub\psi(E^{\bot}), \forall E^{\bot}\in\lat\a^{*}\}.
\end{align*}
$(\u^{l}_{\phi})^{*}$ is a reflexive right $\a^{*}$-submodule determined
by $\psi$. From (1), it follows that $y\otimes x\in(\u^{l}_{\phi})^{*}$
if and only if there exists $E^{\bot}\in \lat\a^{*}$ such that $y\in
E^{\bot}$ and $x\in\psi_{\sim}(E^{\bot})^{\bot}$. Now we compute
$\psi_{\sim}(E^{\bot})^{\bot}$. It follows from the definition of $\pis$
that
\begin{align*}
\hskip -1.25pc \pis(E^{\bot})^{\bot} &= (\vee\{F^{\bot}\in\lat\a^{*}: \psi(F^{\bot})
\not\geq E^{\bot}\})^{\bot}\\
\hskip -1.25pc &= \wedge\{F\in\lat\a: \phi(F)^{\bot}\not\geq E^{\bot}\}\\
\hskip -1.25pc &= \wedge\{F\in\lat\a: \phi(F)\not\leq E\}.
\end{align*}
\end{enumerate}\vspace{-1.5pc}

\hfill \ab
\end{proof}

\section{Rank-one operators}

In this section, we only consider the reflexive right $\a$-submodule,
and omit the superscript and subscript $r$ in the corresponding
notation. The corresponding results for $*$-reflexive left
$\a$-submodule hold naturally, we leave the details for the interested
readers.

\setcounter{theo}{0}
\begin{theor}[\!] Suppose that $\u_{\phi}$ is a reflexive right
$\a$-submodule determined by $\phi$ in $\homp$ and $\r_{\phi}$ the
rank-one submodule generated by rank-one operators in $\u_{\phi}$. Then
$K\in\lat\ri$ if and only if there exists $E\in\lat\a$ such that $E\leq
K\leq\phi_{*}(E)${\rm ,} where $\phi_{*}(E)=\wedge\{\ps(F): F\in\lat\a,
F\not\leq E\}$.
\end{theor}

\begin{proof} Suppose that $K\in\lat\ri$. Let $E=\vee\{F\in\lat\a: F\leq
K\}$. Then $E\in\lat\a$ and $E\leq K$. Let $F\in\lat\a$ with $F\not\leq
E$. We will show that $K\leq\ps(F)$. Let $y$ be any element of $K$. Now
$F\not\leq K$. So we can choose a vector $e\in F$ and $e\not\in K$.
Since $K\in\lat\ri$, for every vector $f\in\ps(F)^{\bot}$, we have
$(e\otimes f)y=(y, f)e\in K$. But since $e\not\in K$ it follows that
$(y, f)=0$ and $y\in\ps(F)$. Thus $K\leq\ps(F)$ and so
$K\leq\phi_{*}(E)$.

Now suppose that there is a subspace $E\in\lat\a$ with $E\leq
K\leq\phi_{*}(E)$. Let $e\otimes f\in\ri$. By Lemma~2.6(1) there is an
element $F\in\lat\a$ such that $e\in F$ and $f\in\ps(F)^{\bot}$. If
$F\leq E$ then $(e\otimes f)K \sub F\sub E\sub K$. If $F\not\leq E$ then
$K\leq\phi_{*}(E)\leq\ps(F)$ and $(e\otimes f)K=(0)\sub K$. In either
case $K$ is invariant under $e\otimes f$ and $K\in\lat\ri$.\hfill \ab
\end{proof}

Suppose that $\u$ is a reflexive right $\a$-module and $\r$ is the
rank-one submodule of $\u$, it follows from Theorems~2.1 and 3.1 that
$K\in\lat\r$ if and only if for some $E\in\lat\a$, $E\leq
K\leq\tau_{*}(E)$, where $\tau(E)=[\u E]$.

\begin{coro}$\left.\right.$\vspace{.5pc}

\noindent If $\l$ is a subspace lattice{\rm ,} $\r$ is the rank-one
subalgebra of $\alg\l$. Then $K\in\lat\r$ if and only if for some
$E\in\lat\alg\l${\rm ,} $E\leq K\leq E_{*}${\rm ,} where
$E_{*}=\wedge\{F_{-}: F\in\lat\alg\l, F\not\leq E\}$ and
$F_{-}=\vee\{G\in\lat\alg\l: G\not\geq F\}$.
\end{coro}

\begin{proof} In this case, $\tau(E)=[(\alg\l) E]=E$ and
$\ts(E)=\vee\{G\in\lat\alg\l: \tau(G)\not\geq E\}=\vee\{G\in\lat\alg\l:
G\not\geq E\}=E_{-}$ and $\tau_{*}(E)=\wedge\{\ts(F): F\in\lat\alg\l,
F\not\leq E\} =\wedge\{F_{-}: F\in\lat\alg\l, F\not\leq E\}=E_{*}$. The
corollary follows from Theorem~3.1.

\hfill \ab
\end{proof}

If $\l$ is a subspace lattice and $E\in\l$, we define
$E_{-}^{\l}=\vee\{F\in\l: F\not\geq E\}$ and
$E_{*}^{\l}=\wedge\{F_{-}^{\l}: F\in\l, F\not\leq E\}$. The following
proposition is due to Longstaff. It gives a similar characterization of
$\lat\r$ by means of the elements in $\l$.

The next two theorems comprise some of the main results of this paper.

\begin{theor}[\!] Suppose that $\l$ is a subspace lattice and $\r$ is
the rank-one subalgebra of $\alg\l$. The following statements are
equivalent.
\begin{enumerate}
\renewcommand\labelenumi{\rm (\arabic{enumi})}
\item $\alg\l=\alg\lat\r${\rm ;}

\item $\lat\alg\l=\lat\r${\rm ;}

\item $[E, E_{*}]\sub\lat\alg\l$ for any $E\in\lat\alg\l${\rm ,} where
$[E, E_{*}]=\{K\in\p: E\leq K\leq E_{*}\}$.
\end{enumerate}
\end{theor}

\begin{proof}
It is clear that (1) is equivalent to (2), we only need to show that (2)
is equivalent to (3).

\noindent (2)\,$\Rightarrow$\,(3).\ \ By definition, $E_{*}=\wedge\{F_{-}:
F\in\lat\alg\l, F\not\leq E\}$. It follows from the definition of
$F_{-}$ that $E\leq F_{-}$ for any $F\not\leq E$. Thus $E\leq E_{*}$ for
any $E\in\lat\alg\l$, so the symbol $[E, E_{*}]$ is meaningful. For any
$E\in\lat\alg\l$ and $K\in [E, E_{*}]$, it follows from Corollary~3.2
that $K\in\lat\r=\lat\alg\l$. Hence $[E, E_{*}]\sub\lat\alg\l$ for any
$E\in\lat\alg\l$.

\noindent (3)\,$\Rightarrow$\,(2).\ \ Since $\r\sub\alg\l$,
$\lat\alg\l\sub\lat\r$. For any $K\in\lat\r$, it follows from Corollary
3.2 that there is an element $E\in\lat\alg\l$ such that $E\leq K\leq
E_{*}$. So $K\in [E, E_{*}]\sub \lat\alg\l$. Thus $\lat\r\sub\lat\alg\l$
and $\lat\r=\lat\alg\l$.\hfill \ab
\end{proof}

\begin{propo}\hskip -.5pc {\rm (\cite{8}, Proposition~3.2)}$\left.\right.$\vspace{.5pc}

\noindent Suppose that $\r$ is the rank-one subalgebra of $\alg\l$. Then
the subspace $K$ belongs to $\lat\r$ if and only if there is a subspace
$E$ of $\l$ such that $E\leq K\leq E_{*}^{\l}$.
\end{propo}

\begin{coro}$\left.\right.$\vspace{.5pc}

\noindent Suppose that $\l$ is a subspace lattice and $\r$ is the
rank-one subalgebra of $\alg\l$. Then $\l=\lat\r$ if and only if $[E,
E_{*}^{\l}]\sub\l$ for any $E\in\l$.
\end{coro}

\begin{proof} Suppose that $\l=\lat\r$. For any $E\in\l$, we can show
$E\leq E_{*}^{\l}$ similarly as in Theorem~3.3. For $K\in [E,
E_{*}^{\l}]$, it follows from Proposition~3.4 that $K\in\lat\r=\l$. So
$[E, E_{*}^{\l}]\sub\l$ for any $E\in\l$.

Conversely, if $[E, E_{*}^{\l}]\sub\l$ for any $E\in\l$. For any
$K\in\lat\r$, it follows from Proposition~3.4 that there is an element
$E\in\l$ such that $E\leq K\leq E_{*}^{\l}$. So $K\in [E,
E_{*}^{\l}]\sub\l$. Thus $\lat\r\sub\l$. Combining with the fact that
$\l\sub\lat\alg\l\sub\lat\r$, we obtain $\l=\lat\r$.\hfill \ab
\end{proof}

\begin{coro}$\left.\right.$\vspace{.5pc}

\noindent Suppose that $\l$ is a subspace lattice. If $[E, E_{*}]\sub\l$
for any $E\in\l${\rm ,} then $\l$ is a reflexive subspace lattice.
\end{coro}

\begin{proof} From Corollary~3.5, it follows that $\l=\lat\r$. Since
$\l\sub\lat\alg\l\sub\lat\r$, so $\l=\lat\alg\l$ and $\l$ is
reflexive.\hfill \ab
\end{proof}

Proposition~3.4, and its Corollary~3.5, answer the question of which
subspace lattices $\l$ are determined by the rank-one subalgebra of
$\alg\l$ in the sense that $\l=\lat\r$. This proposition was used as the
basis of an abstract, lattice-theoretic, way of constructing reflexive
lattices in \cite{10}. Theorem~3.3 gives a sufficient and necessary
condition for which reflexive algebra $\alg\l$ is determined by the
rank-one subalgebra of $\alg\l$ in the sense that $\alg\l=\alg\lat\r$.
In the following, we will consider another sense that $\alg\l$ is
determined by the rank-one subalgebra of $\alg\l$.

\begin{theor}[\!] Suppose that $\a$ is a unital algebra in $\bh$ and
$\u$ is a reflexive right $\a$-submodule{\rm ,} and that $\r$ the rank-one
submodule of $\u$. Then $\ref\r=\u$ if and only if $\tau=\tss${\rm ,} where
$\tau(E)=[\u E]$ for any $E\in\lat\a$.
\end{theor}

\begin{proof}{\it Necessity.}\ \ Recall that the reflexive hull
$\ref\r=\{T\in\bh: Tx\in [\r x], \forall x\in\h\}$ for any $E\in\lat\a$
and $e\otimes f\in\u$. We first show
\begin{equation*}
(e\otimes f)E\sub\tss(E)=\vee\{F\in\lat\a: \ts(F)\not\geq E\}.
\end{equation*}
By virtue of Lemma~2.6(1), there is an element $L\in\lat\a$ such that
$e\in L$ and $f\in\ts(L)^{\bot}$. If $\ts(L)\geq E$ then
\begin{equation*}
(e\otimes f)E=L(e\otimes f)\ts(L)^{\bot}E=(0)\sub\tss(E);
\end{equation*}
if $\ts(L)\not\geq E$ then $L\leq\tss(E)$. Thus
\begin{equation*}
(e\otimes f)E=L(e\otimes f)\ts(L)^{\bot}E\sub L\sub\tss(E).
\end{equation*}
So each rank-one operator of $\u$ maps $E$ into $\tss(E)$ for any
$E\in\lat\a$. For any $A\in\u$ and $x\in E (\in\lat\a)$, since
$\u=\ref\r=\{T\in\bh: Tx\in [\r x], \forall x\in\h\}$, so $Ax\in [\r
x]\sub [\r E]\sub\tss(E)$ and $AE\sub\tss(E)$. Accordingly, $\tau(E)=[\u
E]\sub\tss(E)$ and $\tau\leq\tss$. For any $E\in\lat\a$, it follows from
the definitions that
\begin{equation}
\tss(E)=\vee\{G\in\lat\a: \ts(G)\not\geq E\}
\end{equation}
and
\begin{equation}
\ts(G)=\vee\{F\in\lat\a: \tau(F)\not\geq G\}.
\end{equation}
For $G\in\lat\a$ and $\ts(G)\not\geq E$, if $\tau(E)\not\geq G$ then it
follows from (2) that $\ts(G)\geq E$. This contradiction shows that
$\tau(E)\geq G$. Thus eq.~(1) tells us that $\tau\geq\tss$. Hence
$\tau=\tss$.

{\it Sufficiency.}\ \ Suppose that $\tau=\tss$. It is clear that
$\ref\r\sub\u$, it suffices to show that $\ref\r\supseteq\u$. Suppose
that $A\in\u$ and $x\in\h$. From the definition $\ref\r$, we only need
to prove that $Ax\in [\r x]$.

Define $E$ by $E=\wedge\{F\in\lat\a: x\in F\}$. Observe that the
intersection is over a non-empty family of subspaces of $\lat\a$ since
$x\in\h$. Clearly $x\in E$ and $E\in\lat\a$. By the hypothesis, 
\begin{equation*}
\tau(E)=[\u E]=\vee\{G\in\lat\a: \ts(G)\not\geq E\}
\end{equation*}
and hence the set of all $G\in\lat\a$ with $\ts(G)\not\geq E$ has a
dense linear span in $[\u E]$. Therefore for any $\epsilon >0$, there is
a finite set $G_{i} (1\leq i\leq n)$ of subspaces of $\lat\a$ with
$\ts(G_{i})\not\geq E$ and a set of vectors $x_{i}\in G_{i} (1\leq i\leq
n)$ with the property that
\begin{equation*}
\|Ax-(x_{1}+\cdots +x_{n})\|<\epsilon.
\end{equation*}

The definition of $E$ and the condition $\ts(G_{i})\not\geq E (1\leq
i\leq n)$ implies that $x\not\in\ts(G_{i}) (1\leq i\leq n)$ and so there
exists $y_{i} \in\ts(G_{i})^{\bot}$ with
\begin{equation*}
\langle x, y_{i} \rangle \not=0,\qq\forall 1\leq i\leq n.
\end{equation*}
By suitably scaling $y_{i}$ if needed we may assume that $\langle x,
y_{i} \rangle =1$ and so $(x_{i}\otimes y_{i})x=x_{i}$ for $1\leq i\leq
n$. By Lemma~2.6(1), $x_{i}\otimes y_{i}\in\r$. Thus
\begin{equation*}
\left\Vert Ax- \left( \sum\limits_{i=1}^{n}x_{i}\otimes y_{i} \right)
x\right\Vert = \|Ax-(x_{1}+\cdots+x_{n})\|<\epsilon,
\end{equation*}
and this shows that $Ax\in [\r x]$ and $A\in\ref\r$. Hence
$\ref\r=\u$.\hfill \ab
\end{proof}

\begin{coro}$\left.\right.$\vspace{.5pc}

\noindent Suppose that $\l$ is a subspace lattice and $\r$ is the
rank-one subalgebra of $\alg\l$. Then $\l$ is completely distributive if
and only if $\ref\r=\alg\l$.
\end{coro}

\begin{proof} In this case, $\tau(E)=[(\alg\l)E]=E$, $\ts(E)=E_{-}$ and
$\tss(E)=\vee\{F\in\lat\alg\l: \ts(F)\not\geq
E\}=\vee\{F\in\lat\alg\l: F_{-}\not\geq E\}=E_{\sharp}$. Now
suppose that $\ref\r=\alg\l$. So from Theorem~3.7, it follows that
$E=E_{\sharp}$ for any $E\in\lat\alg\l$. Theorem~5.2 of ref.~\cite{8} shows
that $\lat\alg\l$ is completely distributive. $\l \sub\lat\alg\l$
implies that $\l$ is also completely distributive.

Conversely, suppose that $\l$ is completely distributive. From
(\cite{8}, Theorem~6.1), it follows that $\l=\lat\alg\l$. So
$E=E_{\sharp}$ for any $E\in\l=\lat\alg\l$. It follows from Theorem~3.7
that $\ref\r=\alg\l$.\hfill \ab
\end{proof}

Corollary~3.8 was first proved by Lambrou~(\cite{6}, Theorem~3.1). From
the above proof, we can easily obtain that $\ref\r=\alg\l$ is equivalent
to the complete distributivity of $\lat\alg\l$. Thus it follows from
(\cite{8}, Theorem~5.2) that $\ref\r=\alg\l$ if and only if $E=E_{*}$
for any $E\in\lat\alg\l$. Comparing with Theorem~3.3, shows the
differences between $\alg\lat\r$ and $\ref\r$.

\end{document}